\documentclass{amsart}
\usepackage{eufrak}
\usepackage{amssymb,amsmath}
\author{mihai popa}
\address{UCLA, Los Angeles, CA 90095-1555}
\email{mvpopa@math.ucla.edu}
\title{A Non-commutative Analogue of Gaussian Hilbert Spaces}
\newtheorem{defin}{Definition}[section]
\newtheorem{claim}{}[section]
\newtheorem{thm}[claim]{Theorem}
\newtheorem{remark}{Remark}[section]
\newtheorem{prop}[claim]{Proposition}
\newtheorem{cor}[claim]{Corollary}
\newcommand{\cA}{\mathcal{A}}

\begin{document}
\begin{abstract}
 The paper gives analogues of some starting results in the theory of Gaussian Hilbert Spaces for semicircular distributed random variables. The transition from the commutative to the free frame is done considering matrices of increasing dimension and utilizing the Amitsur-Levitzki Theorem.
\end{abstract}
\maketitle
\section{Introduction}
  The Gaussian Hilbert spaces (i.e. complete linear spaces of random variables with centered Gaussian distributions) are referred in many works (see \cite{Jan}) as an important structure in probability theory, stochastic processes and fields, stochastic integration, quantum field theory etc. Two well-know basic results in this topic are that a vector is Gaussian if and only if any linear combination of its components is also Gaussian and that any Gaussian Hilbert space can be seen as the space of sums of creation and annihilation operators on the Symmetric Fock Space over some unique, up to an isometry, Hilbert space.

 This paper gives analogues of the results above for semicircular distributed random variables which play in the Voiculescu's Free Probability theory a similar role of Gaussian variables in the classical probability frame. The main idea in the proof of the classical results is identifying the scalar coefficients of some polynomials in several commuting variables. In the non-commutative frame matrices of increasing dimension are replacing the scalars and while the technique deployed many times in Free Probability theory for obtaining analogous results to classical statements uses the property of matrices with Gaussian distributed random entries and increasing dimension to asymptotically behave as free semicircular variables, here the transition from commutative to free frame is done via the Amitsur-Levitzki Theorem. 

  Sections 2. and 3. contain preliminary definitions and constructions, primarily recalling the ones from \cite{Sp2} and \cite{Pi}; the main results are given in Section 3.

\section{B-probability spaces and B-semicircular random variables}

\begin{defin}
 Let B be a C*-algebra. By a B-probability space, or a  non-commutative probability space over B, we will understand a pair $(\cA, \varphi)$ consisting of an algebra $\cA$ over B and a conditional expectation 
$ \varphi:\cA\longrightarrow B $. 
If $\cA$ is a $\ast$-algebra, the pair $(\cA, \varphi)$ will be called $\ast$-B-probability space.
\end{defin}
Consider the moment function 
$\Hat{\mu}=(\mu^{(n)})_{n\in\mathbb{N}}$, 
\[\mu^{(n)}:\underbrace
{\mathcal{A}\otimes_{B}\cA\otimes_{B}\dots\otimes_{B}\mathcal{A}}
_{\text{n  times}}\longrightarrow B\]
given by 
$\mu^{(n)}({a_1}\otimes\dots\otimes{a_n})=\varphi({a_1}\dots{a_n})$ 
 and the correspondent cumulant function 
$\Hat{k}=(k^{(n)})_{n\in\mathbb{N}}$ 
\[k^{(n)}:\underbrace{\cA\otimes_{B}\mathcal {A}\otimes_{B}\dots\otimes_{B}\mathcal{A}}
_{\text{n  times}}\longrightarrow B\]
given by 
 \begin{equation}
\label{eq:1}
\mu^{(n)}({a_1}\otimes\dots\otimes{a_n})=\sum_{\pi\in NC(n)} {k_\pi}[{a_1}\otimes\dots\otimes{a_n}]
\end{equation}
 where $NC(n)$ is the set of all non-crossing partitions on $\{1,2,\dots,n\}$ and \\${k_\pi}[{a_1}\otimes\dots\otimes{a_n}]$ is defined recursively by:

\begin{enumerate}

\item[(a)]
$k_{\pi_{1}\cup{\pi_{2}}}={k_{\pi_{1}}}{k_{\pi_{2}}}
\text{, for any }\pi_{1}\in{NC(1,2,\dots,n-j)}, \pi_{2}\in{NC({n-j+1},\dots,{n})}$
\medskip
\item[(b)]
$k_{1\negthickspace{1}}[{a_{1}}\otimes\dots\otimes{a_{n}}]=k^{(n)}({a_{1}}\otimes\dots\otimes{a_{n}})$ , for  ${1\negthickspace{1}}$ the partition with a single equivalence class 
\medskip
\item[(c)]
$k_{\pi}[{a_{1}}\otimes\dots\otimes{a_{p}}\otimes{a_{p+1}}\otimes\dots\otimes{a_{p+r}}
\otimes\dots\otimes{a_{n}}]$\
$=k^{(n_r)}({a_{1}}\otimes\dots\otimes{a_{p}}\otimes{k_{\sigma}
[{a_{p+1}}\otimes\dots\otimes{a_{p+r}}]}{a_{p+r+1}}\otimes\dots\otimes{a_{n}}),$ 
for any 
${\sigma}\in {NC(p+1,\dots,p+r)}$ 
and $\pi$ the partition obtaining by 
${\sigma}\cup{(a_{1},\dots a_{p},a_{p+r},\dots{a_{n}})}$

\end{enumerate}

For illustrations, the explicit formulas for connections between the small order moment and cumulant functions
are

\begin{eqnarray*}
\mu^{(1)}(a_{1})&=&k^{(1)}(a_1)=\varphi(a_{1})\\
& &\\
\mu^{(2)}(a_{1}\otimes{a_{2}})
&=&
k_{\{(1,2)\}}[a_{1}\otimes{a_{2}}]+k_{\{(1),(2)\}}[a_{1}\otimes{a_{2}}]\\
&=&
k^{(2)}(a_{1}\otimes{a_{2}})+k^{(1)}(a_{1})k^{(1)}(a_{2})\\
& &\\
k^{(2)}(a_{1}\otimes{a_{2}})
&=&
\mu^{(2)}(a_{1}\otimes{a_{2}})-\mu^{(1)}(a_{1})\mu^{(1)}(a_{2})\\
&=&
\varphi(a_{1}a_{2})-\varphi(a_{1})\varphi(a_{2})\\
& &\\
\mu^{(3)}(a_{1}\otimes{a_{2}}\otimes{a_{3}})
&=&
k_{\{(1,2,3)\}}[a_{1}\otimes{a_{2}}\otimes{a_{3}}]+k_{\{(1),(2,3)\}}
[a_{1}\otimes{a_{2}}\otimes{a_{3}}]\\
 & & +k_{\{(1,2),(3)\}}[a_{1}\otimes{a_{2}}\otimes{a_{3}}]+k_{(1,3),(2)}
[a_{1}\otimes{a_{2}}\otimes{a_{3}}]\\ 
& & +k_{\{(1),(2),(3)\}}[a_{1}\otimes{a_{2}}\otimes{a_{3}}]\\
&=&
k^{(3)}(a_{1}\otimes{a_{2}}\otimes{a_{3}})+k^{(1)}(a_{1})k^{(2)}({a_{2}}
\otimes{a_{3}})\\ 
& & +k^{(2)}({a_{1}}\otimes{a_{2}})k^{(1)}(a_{3})+k^{(2)}({a_{1}}
\otimes{k^{(1)}(a_{2})}{a_{3}})\\
& & +k^{(1)}(a_{1})k^{(1)}(a_{2})k^{(1)}(a_{3})\\
& &\\
k^{(3)}(a_{1}\otimes{a_{2}}\otimes{a_{3}})
&=&
\mu^{(3)}(a_{1}\otimes{a_{2}}\otimes{a_{3}})-
\mu^{(1)}(a_{1})\mu^{(2)}({a_{2}}\otimes{a_{3}})\\
& & -\mu^{(2)}({a_{1}}\otimes{\mu^{(1)}(a_{2})}{a_{1}})-\mu^{(2)}({a_{1}}\otimes{a_{2}}\mu^{(1)}(a_{3})\\
& & +2\mu^{(1)}(a_{1})\mu^{(1)}(a_{2})\mu^{(1)}(a_{3})\\
& &=\varphi(a_{1}{a_{2}}{a_{3}})-\varphi(a_{1})\varphi({a_{2}}{a_{3}})-\varphi({a_{1}}{\varphi(a_{2})}{a_{3}})\\
& & -\varphi({a_{1}}{a_{2}})\varphi(a_{3})+ 2\varphi(a_{1})\varphi(a_{2})\varphi(a_{3})
\end{eqnarray*}
\medskip
 
Given a family $(X_{1},\dots,X_{m})$ of elements from $\cA$, the information about their joint moments with respect to $\varphi$ is encoded in the collection of mappings   \[{\xi}={({\xi}_{n;i_{0},\dots,i_{n}})_{n\in\mathbb{N};
i_{0},\dots,i_{n}\in\lbrace 1,\dots,m\rbrace}},\] 
where 
${\xi}_{n;i_{0},\dots,i_{n}}:\underbrace{{B}\times\dots\times{B}}\longrightarrow{B}$ 
are given by

\begin{eqnarray}
{\xi}_{n;i_{0},\dots,i_{n}}(b_{1},\dots,b_{n})
&=&
{{k}^{(n+1)}({X}_{i_{0}}\otimes{b}_{1}{X}_{i_{1}}\otimes\dots\otimes{b}_{n}{X}_{i_{n}})}
\nonumber\\
&=&
{{k}^{(n+1)}({X}_{i_{0}}{b}_{1}\otimes{X}_{i_{1}}{b}_{2}\otimes\dots\otimes{X}_{i_{n}})}
\nonumber
\end{eqnarray}
or, equivalently, in the collection 
$${\eta}={({\eta}_{n;i_{0},\dots,i_{n}})_{n\in\mathbb{N};
i_{0},\dots,i_{n}\in\lbrace 1,\dots,m\rbrace}}$$
 where
\begin{eqnarray}
{\eta}_{n;i_{0},\dots,i_{n}}(b_{1},\dots,b_{n})&=&{{\mu}^{(n+1)}({X}_{i_{0}}\otimes{b}_{1}{X}_{i_{1}}\otimes\dots\otimes{b}_{n}{X}_{i_{n}})}\nonumber\\
&=&{{\mu}^{(n+1)}({X}_{i_{0}}{b}_{1}\otimes{X}_{i_{1}}{b}_{2}\otimes\dots\otimes{X}_{i_{n}})}\nonumber
 \end{eqnarray}
since moments are polynomials in cumulants and viceversa.

 If $X_{1}=X_{2}=\dots=X_{m}=X$, we will write ${\xi}_{n;X}$, respectivelly ${\eta}_{n;X}$ for the functionals above.

\begin{remark}
In the case $B=\mathbb{C}$, one only needs the coeficients ${\xi}_{n;i_{0},\dots,i_{n}}(1,\dots,1)$ 
(respectively the coeficients 
${\eta}_{n;i_{0},\dots,i_{n}}(1,\dots,1)$
), since 
$${\xi}_{n;i_{0},\dots,i_{n}}(b_{1},\dots,b_{n})
={b}_{1}\dots{b}_{n}{\xi}_{n;i_{0},\dots,i_{n}}({1},\dots,{1})$$
 (and the analogous relation for $\eta$'$s$).
\end{remark} 

\begin{defin}
 For $(\cA, \varphi)$ a $\ast$-B-probability space, a selfadjoint element $x$ from $\mathcal{A}$ is called B-semicircular if ${\xi}_{n;x}=0$ for all $n\neq{1}$.
\end{defin}

\begin{remark}
 The relation (\ref{eq:1}) enable us to formulate the above definition in terms of moments: $x={x}^{*}\in\mathcal{A}$ 
is B-semicircular if
\[{\eta_{n;x}(b_{1},\dots,b_{n})}
={\mu}^{(n+1)}(x\otimes{b_1}x\otimes\dots\otimes{b_n}x)
=\sum_{\pi\in NCPP(n)}{{k_{\pi}}[x\otimes{b_1}x\otimes\dots\otimes{b_n}x]}\]
 where NCPP(n) denotes all the pair-partitions from NC(n).
\end{remark}
 \begin{defin}\text{(see \cite{EP})}
 The semicircular elements \hspace{0.5mm} $x_1,\dots,x_m$\  from $(\cA, \varphi)$, a non-commutative probability space,  are said to satisfy the 0-Wick Theorem if all  their coefficients ${\xi}_{p;i_{0},\dots,i_{n}}(1,\dots,1)$ are zero unless $p=1$, or,  equivalently, if their joint moments are computed by the formula:
\begin{equation}
 {\varphi}(x_{i_1}\dots{x_{i_p}})
=\sum_{\gamma\in{NCPP(p)}}(\prod_{(k,l)\in{\gamma}}\varphi(x_{i_k}x_{i_l})),\ 
\text{for any}\  i_1,\dots,i_p\in\{1,\dots,m\}
\nonumber
\end{equation}
\end{defin}

\section{On a Particular Hilbert Bimodule }

Let $\mathcal{H}^0$ be a real Hilbert space and $\mathcal{H}$ its complexification. For $f\in\mathcal{H}^0$ we denote by 
$a_{1}(f)$,\ $a_{1}^{*}(f)$,
\  $G_{1}(f)$,\ $a_{0}(f)$,\ $a_{0}^{*}(f)$,
\  $G_{0}(f)$ 
the CCR operators of annihilation, creation and Gaussian (sum of the creation and annihilation operators) over $f$, respectively their full Fock Space analogues. Let $\Gamma_j$ be the von Neumann Algebra generated by
 $\{G_j(f),\ f\in\mathcal{H}\}$, $j=0,1$. 
In both situations, the vacuum state will be denoted by 
$\langle{\Omega},\cdot{\Omega}\rangle$.\\ 
\indent
 Denote by $\mathcal{K}(l^{2}(\mathbb{N}))$ the C*-algebra of compact operators on $l^{2}(\mathbb{N})$; for any $n\in\mathbb{N}$, $M_{n}(\mathbb{C})$ will be seen as a subalgebra of  
$\mathcal{K}(l^{2}(\mathbb{N}))$
 via the identification to the set of elements with non-zero entries only in the upper left $n\times n$-corner\\
\indent
 For 
$(\mathcal{A}, \varphi)$ 
a non-commutative probability space  and $\widetilde{\mathcal{A}}={\mathcal{K}(l^{2}(\mathbb{N}))}\otimes{\mathcal{A}}$, 
(respectively 
$M_{n}(\mathcal{A})={M}_{n}(\mathbb{C})\otimes{\mathcal{A}}$), $\widetilde{\varphi}
=id_{\mathcal{K}(l^{2}(\mathbb{N}))}\otimes{\varphi}$ 
(respectively \ 
$\varphi_{n}=id_{{M}_{n}(\mathbb{C})}\otimes{\varphi}$\ 
the  
${\mathcal{K}(l^{2}(\mathbb{N}))}$, 
respectively 
${M}_{n}(\mathbb{C})$-probability spaces 
$(\widetilde{\mathcal{A}}, \widetilde{\varphi})$ 
and 
$({M}_{n}(\mathcal{A}), \varphi_{n})$ 
will be considered with the bimodule structure given by the linear extension of the action
\[Ax=xA=A\otimes{x},\  
\text{for any $x\in\cA$ and \ $A$ in  
${\mathcal{K}(l^{2}(\mathbb{N}))}$
 or in 
${{M}_{n}(\mathbb{C})}$}\].
\begin{defin}
  For $\mathcal{B}$ a C*-algebra and $\mathcal{E}$ a Hilbert $\mathcal{B}$-bimodule, the Full Fock Bimodule  over $\mathcal{E}$ is the semi-inner-product $\mathcal{B}$-bimodule  
\[\mathcal{F}(\mathcal{E})
=\mathcal{B}\widetilde{\Omega}\oplus\mathcal{E}\oplus(\mathcal{E}
\otimes_{\mathcal{B}}\mathcal{E})\oplus(\mathcal{E}
\otimes_{\mathcal{B}}\mathcal{E}\otimes_{\mathcal{B}}\mathcal{E})\oplus\dots\]
equipped with the pairing
\[\langle\cdot,\cdot\rangle{:}\mathcal{F}(\mathcal{E})\times\mathcal{F}(\mathcal{E})
\rightarrow{\mathcal{B}}\]
given by the linear extension of (
$A,B\in{\mathcal{B}}, \widetilde{x_j},\widetilde{y_j}\in\mathcal{E}$):

\begin{eqnarray*}
\langle{A},{B}\rangle
&=&
A^{*}B\\
\langle{\widetilde{x_1}\otimes\dots\widetilde{x_k}},
{\widetilde{y_1}\otimes\dots\widetilde{y_l}}\rangle
&=&
\delta_{k}^{l}
\langle{\widetilde{x_k}}\langle\widetilde{x_1}
\otimes\dots\widetilde{x_{k-1}},
\widetilde{y_1}\otimes\dots\widetilde{y_{k-1}}\rangle
{\widetilde{y_k}}\rangle\\
&=&
\delta_{k}^{l}
\langle{\widetilde{{x_k}}},
\langle{\widetilde{{x_{k1}}},\dots,\langle{\widetilde{x_1}},\widetilde{{y_1}}}
\rangle\dots\widetilde{y_{k-1}}\rangle\widetilde{y_k}\rangle
\end{eqnarray*}
The Symmetric Fock Bimodule $\mathcal{S}(\mathcal{E})$ over $\mathcal{E}$ is the submodule of $\mathcal{F}(\mathcal{E})$ generated by the elements of the form $\underbrace{\widetilde{a}\otimes\widetilde{a}\otimes\dots\otimes\widetilde{a}}_{n\ times}$\ ,
\ $a\in\mathcal{E}$, $n\in\mathbb{N}$.
\end{defin}

  For\hspace{5pt}  
$\widetilde{x}\in\mathcal{E}$ 
we consider, as in \cite{Pi}, 
$G(\widetilde{x})=a(\widetilde{x})+a^{*}(\widetilde{x})$
, where 
$a(\widetilde{x})$ and $a^{*}(\widetilde{x})$ 
are respectively the annihilation and the creation operators on  
$\mathcal{F}(\mathcal{E})$ 
given by 
\begin{enumerate}
\item[]
$a(\widetilde{x})A=0$ , $A\in{\mathcal{B}}$
\item[]
$a(\widetilde{x})\widetilde{x_1}\otimes\dots\widetilde{x_m}
=\langle\widetilde{x},\widetilde{x_1}\rangle\widetilde{x_2}\otimes\dots\widetilde{x_m}$ , 
$\widetilde{x_j}\in{\mathcal{E}}$
\end{enumerate}

\begin{enumerate}
\item[]
$a^{*}(\widetilde{x})A=\widetilde{x}A$ , $A\in{\mathcal{B}}$
\item[]
$a^{*}(\widetilde{x})\widetilde{x_1}\otimes\dots\widetilde{x_m}
=\widetilde{x}\otimes\widetilde{x_1}\otimes\dots\widetilde{x_m}$ , 
$\widetilde{x_j}\in{\mathcal{E}}$
\end{enumerate}\
 
   Denote by $\widetilde{H}$ the Hilbert 
$\mathcal{K}(l^2(\mathbb{N}))$-bimodule 
$\mathcal{K}(l^2(\mathbb{N}))\otimes{H}$
, where 
$Ax=xA=\\ 
A\otimes{x},\ \text{for any}\ x\in{H}\ \text{and}\ A\in\mathcal{K}(l^2(\mathbb{N}))$
, and denote by 
$\widetilde{H}^0$ 
be the set of all selfadjoint elements form $\widetilde{H}$ 
with respect to  the involution given by 
$A\otimes{x}=Ax\mapsto{A^*{x}}={A^*\otimes\Bar{x}}$. 
Let $\mathfrak{B}$ be the algebra generated by 
$\mathcal{K}(l^{2}(\mathbb{N}))$  \ and 
$\{G(\widetilde{x}) , \widetilde{x}\in\widetilde{H^0}\}$. 
Then $\mathfrak{B}$ with the 
$\mathcal{K}(l^{2}(\mathbb{N}))$-valued functional $\langle\widetilde{\Omega},\cdot\widetilde{\Omega}\rangle$ 
is an 
$\mathcal{K}(l^{2}(\mathbb{N}))$-$\ast$-probability space  and we have that:
\begin{prop}
 For any $\widetilde{x}\in{\widetilde{H^0}}$, the operator $G(\widetilde{x})$ is $\mathcal{K}(l^{2}(\mathbb{N}))$-semicircular.
\end{prop}
\begin{proof}
 It is enough to prove the assertion for elements of the form  $\widetilde{x}=\sum_{j=1}^{m}A_{j}{x_j}$
 where 
$x_j\in{\mathcal{H}^0}$, $n\in\mathbb{N}$
 and 
$A_j=A_j^{*}\in{M_{n}(\mathbb{C})}$
. Then 
$G(\widetilde{x})=\sum_{j=1}^{m}A_{j}G_0({x_j})$, 
where  $A_{j}G_0({x_j})$ 
denotes the operator $A_{j}\otimes{G_0({x_j})}$\ acting on $\mathcal{K}(l^2(\mathbb{N})\otimes{\mathcal{F}(H)}$\ 
and for any 
$B_j\in\mathcal{K}(l^{2}(\mathbb{N}))$ 
one has that
 
\begin{eqnarray*}
k^{(p+1)}(G(\widetilde{x})\otimes& &
\mspace{-46mu}
{B_1}G(\widetilde{x})\dots\otimes{B_p}G(\widetilde{x}))=\\
&=&
k^{(p+1)}(\sum_{j=1}^{m}A_{j}G({x_j})\otimes{B_1}\sum_{j=1}^{m}A_{j}G({x_j})
\dots\otimes{B_p}\sum_{j=1}^{m}A_{j}G({x_j}))\\
&=&
\sum_{\substack{(j_0,\dots{j_p})\\1\leq{j_k}\leq{m}}}A_{j_0}{B_1}{A_{j_1}}
\dots{B_p}{A_{j_p}}k^{(p+1)}(G_0(x_{j_o})\otimes
\dots\otimes{G_0({x}_{j_p})})=0 
\end{eqnarray*}
  unless $p=1$, since, from [1], 
$G_0(x_{j_o}),\dots,{G_0({x}_{j_p})}$ 
form a semicircular family.
\end{proof}
\begin{remark}
 For 
$(H_{j}, \xi_{j})_{j\in\mathbb{N}}$ 
a family of real Hilbert spaces with distinguished unit vectors $\xi_{j}$, one has:
\begin{equation*}
 \mathcal{F}(\widetilde{H_j})\cong\mathcal{S}(\widetilde{H_j})
\end{equation*}
and
\begin{equation*}
 \underset{j\in\mathbb{N}}{\text{\Large$\ast$}}
\mathcal{F}(\widetilde{H_j})\cong\bigotimes_{j\in\mathbb{N}}\mathcal{F}(\widetilde{H_j})
\end{equation*}
where {\text{\Large$\ast$}} denotes the free product of Hilbert bimodules, (as considered in \cite{VDN} for Hilbert spaces) with amalgamation over $\mathcal{K}(l^{2}(\mathbb{N}))$
\end{remark}
\begin{proof}
 Since 
$\mathcal{F}(\underset{j\in\mathbb{N}}{\bigoplus}\widetilde{H_j})
={\underset{j\in\mathbb{N}}{\text{\Large$\ast$}}}\mathcal{F}(\widetilde{H_j})$ 
and $\mathcal{S}(\underset{j\in\mathbb{N}}{\bigoplus}\widetilde{H_j})
=\underset{j\in\mathbb{N}}{\bigotimes}\mathcal{S}(\widetilde{H_j})$ 
(same arguments as in \cite{VDN} and \cite{Gui}), it suffices to prove that $\mathcal{F}(\widetilde{H})=\mathcal{S}(\widetilde{H})$ 
for $H$ a real Hilbert Space, i.e. 
$\mathcal{S}(\widetilde{H})$ 
contains elements of the form 
$Aa_1\otimes{a}_2\otimes\dots\otimes{a}_n$ 
for arbitrary 
$A\in {M_m({\mathbb{C}})}\subset{\mathcal{K}(l^{2}(\mathbb{N}))}$
 and 
$a_1,\dots, a_n\in{H}$.\\
 \indent Consider 
$\widetilde{a}=AA_1a_1+A_{2}a_{2}+\dots+A_na_n$, where $A_j\in{M}_{n}(M_m(\mathbb{C}))\subset{\mathcal{K}(l^{2}(\mathbb{N}))}$
 are given by
\begin{displaymath}
A_j=(\alpha_{k,l}^{(j)})_{k,j=1}^{n}\ \text{for}\ \alpha_{k,l}^{(j)}=
\begin{cases}
\delta_{j}^{k}\delta_{j+1}^{l}I_m& \text{if}\ j=1,\dots,n-1 
\cr \delta_{n}^{k}\delta_{1}^{l}I_m& \text{if}\ j=n 
\end{cases}
\end{displaymath}

 Then, for $\sigma\in{S}_{n}$ and 
$B_{j}=\begin{cases}AA_1,&\text{if}\ j=1 \cr A_j, &\text{if}\ j=2,\dots,n\end{cases}$ 
one has that\\ 
$B_{\sigma(1)}B_{\sigma(2)}\dots{B}_{\sigma(n)}=\begin{cases}A& 
\text{if}\ \sigma=1_n \cr0&\text{if}\ \sigma \neq 1_n \end{cases}$
\begin{eqnarray*}
 \underbrace{\widetilde{a}\otimes\widetilde{a}\dots
\otimes\widetilde{a}}_{n times}&=&\sum_{\sigma\in{S}_{n}}B_{\sigma(1)}a_{\sigma(1)}
\otimes B_{\sigma(2)}a_{\sigma(2)}\otimes
\dots\otimes B_{\sigma(n)}a_{\sigma(n)}\\
&=&
\sum_{\sigma\in{S}_{n}}B_{\sigma(1)}B_{\sigma(2)}\dots{B}_{\sigma(n)}a_{\sigma(1)}\otimes
\dots{a}_{\sigma(n)}\\
&=&
Aa_1\otimes{a}_2\otimes\dots\otimes{a}_n
\end{eqnarray*}
 \end{proof}
\bigskip

\section{Semicircular and complex semicircular families}

 By a Gaussian Hilbert space it will be understood a closed linear subspace in ${{L}^{2}_{\mathbb{R}}}(\Omega,\Sigma,P)$\ where any element is a centered Gaussian random variable, for $(\Omega,\Sigma,P)$  a (classical) probability space. Is said that the random variables $X_{1},\dots,X_{n}$ are jointly Gaussian or form a Gaussian family if the vector $(X_{1},\dots,X_{n})$ is Gaussian-distributed, or, equivalently, the joint moments of \ $X_{1},\dots,X_{n}$ are computed according to the Wick Theorem (another equivalent condition is that the joint classical cumulants of $X_{1},\dots,X_{n}$ are all zero, except the ones of order 2 (See \cite{Jan})).

 Consider a non-commutative probability space $(\mathcal{A},\varphi)$; the (selfadjoint) elements \  $x_{1},\dots,x_{n}\in\mathcal{A}$ are said to form a semicircular family if the vector $(x_{1},\dots,x_{n})$ is the limit of   a (non-commutative) Free Central Limit Theorem, or equivalently, the joint moments of $x_{1},\dots,x_{n}$ are computed according to the 0-Wick Theorem (see \cite{Sp2}).   

The following result is the analogue of the classical theorem (cf. \cite{Jan}): \\
  \textit{The vector-valued random variable} 
$(X_{1},\dots,X_{n})$ \textit{is Gaussian if and only if for any} 
$\alpha_{j}\in\mathbb{R}$ ($j=1,\dots,n)$,  
\textit{the random variable} 
$\sum_{j=1}^{n}{\alpha_{j}{X}_{j}}$ 
\textit{is Gaussian}.

\begin{thm}\label{thm1}
 Let $(\cA,\varphi)$ be a non-commutative $\ast$-probability space and 
$x_{1},\dots,x_{n}\in\cA$.
 Then ($x_{1},\dots,x_{n}$) form a semicircular family if and only if\  $\widetilde{x}=\sum_{j=1}^{n}{A_{j}x_{j}}$
 is $\mathcal{K}(l^{2}(\mathbb{N}))$-semicircular in 
$(\widetilde{\mathcal{A}}, \widetilde{\varphi})$, 
for any $A_{j}\in \mathcal{K}(l^{2}(\mathbb{N}))_{SA}$. 
\end{thm}

\begin{proof}
Suppose ($x_{1},\dots,x_{n}$) form a semicircular family. Then, for $p\neq1$, one has that:

 \begin{eqnarray}
{\xi}_{p,\widetilde{x}}(B_{1},\dots,B_{p})
&=&
{k^{(p+1)}(\widetilde{x}\otimes{B}_{1}\widetilde{x}\otimes
\dots\otimes{B}_{p}\widetilde{x})}
\nonumber\\
&=&
{k}^{(p+1)}(\sum_{j=1}^{n}{A_{j}{x}_{j}}\otimes{B_1}\sum_{j=1}^{n}{A_{j}{x}_{j}}\otimes
\dots\otimes{B_p}\sum_{j=1}^{n}{A_{j}{x}_{j}})
\nonumber\\
&=&
{k}^{(p+1)}(\sum_{(j_{0},\dots,j_{p})}{A_{j_0}x_{j_0}\otimes{B_1}A_{j_1}x_{j_1}\otimes
\dots\otimes{B}_{p}A_{j_p}x_{j_p}})
\nonumber\\
&=&
{k}^{(p+1)}(\sum_{(j_{0},\dots,j_{p})}{A_{j_0}B_{1}A_{j_1}\dots{B}_{p}A_{j_p}}x_{j_o}
\otimes\dots\otimes{x}_{j_p})
\nonumber\\
&=&
\sum_{\substack{(j_{0},\dots,j_{p})\\1\leq{j_k}\leq{n}}}
{A_{j_0}B_{1}A_{j_1}\dots{B}_{p}A_{j_p}}{k}^{(p+1)}(x_{j_o}\otimes
\dots\otimes{x}_{j_p})=0
\nonumber
\end{eqnarray}
since 
${k}^{(p+1)}(x_{j_o}\otimes\dots\otimes{x}_{j_p})=0$
 for $p\neq1$, so q.e.d..\bigskip

Conversely, take 
$\widetilde{x}=\sum_{j=1}^{n}{A_{j}x_{j}}$
\ for $A_{j}\in{M}_m(\mathbb{C})_{SA}$. If\ $\widetilde{x}$ is 
$\mathcal{K}(l^{2}(\mathbb{N}))$-semicircular, then 
$\xi_{p;\widetilde{x}}(I_m,\dots,I_m)=0$ for $p\neq1$, so
\[0={k}^{(p+1)}(\widetilde{x}\otimes\dots\otimes\widetilde{x})
={k}^{(p+1)}(\sum_{j=1}^{n}{A_{j}{x}_{j}}\otimes
\dots\otimes\sum_{j=1}^{n}{A_{j}{x}_{j}})\]
\[=\sum_{\substack{(j_{0},\dots,j_{p})\\1\leq{j_k}\leq{n}}}
{A_{j_0}\dots{A}_{j_p}}{k}^{(p+1)}(x_{j_0}\otimes
\dots\otimes{x}_{j_p})\]
 The above expression can be seen as a polynomial $f(A_1,\dots,A_p)$ in the non-commutative variables $A_1,\dots,A_p$. Note that $f(A_1,\dots,A_p)=0$ for any selfadjoint $A_1,\dots,A_p\in{M_{m}(\mathbb{C})}\subset\mathcal{K}(l^{2}(\mathbb{N}))$
 and any $m\in\mathbb{N}$, particularly for $m=p+1$.
On the other hand, the Theorem of Amitsur-Levitzki and Lemma I.3.2 from \cite{J} state that the lowest degree of such a polynomial with non-zero coefficients is $2(p+1)$, therefore all the coefficients of $f$ must be zero, i.e. 
${k}^{(p+1)}(x_{j_0}\otimes\dots\otimes{x}_{j_p})=0$ 
for $p\neq1$, so $x_1,\dots,x_n$ form a semicircular family.
\end{proof}
\bigskip
\begin{defin}
 The element $c$ from the non-commutative $\ast$-probability space $(\cA,\varphi)$ will be called \textit{complex semicircular} if there are two jointly semicircular elements 
$s_1, s_2\in\cA$ 
such that $c=s_1+is_2$. If $s_1$ and $s_2$ are free, then $c$ is called circular.
  
The elements $c_{1},\dots,c_{n}\in\mathcal{A}$ are said to form a complex semicircular family if $\Re{c_1}, \Im{c_1},\dots,\Re{c_n}$, $\Im{c_n}$ form a semicircular family. 

 For $(\cA,\varphi)$ a $\ast$-B-probability space, 
$\widetilde{c}\in\cA$ 
is called complex B-semicircular if there are two B-semicircular elements 
$\widetilde{s_1}, \widetilde{s_2}\in\cA$ 
such that 
$\widetilde{c}=\widetilde{s_1}+i\widetilde{s_2}$
 and for any $b_1, b_2\in{B}$ if $b_1\widetilde{s_1}+b_2\widetilde{s_2}$ 
is selfadjoint then it is B-semicircular. 
\end{defin}

 The complex semicircular elements are the free analogue of the classical complex Gaussian random variables.  As in the classical case, $\ref{thm1}$ implies the following:
\begin{cor}
Let $(\cA,\varphi)$ be a non-commutative $\ast$-probability space and  $c_{1},\dots,c_{n}$ be a family of elements in $\cA$.
 Then ($c_{1},\dots,c_{n}$) form a complex semicircular family if and only if\  $\widetilde{x}=\sum_{j=1}^{n}{A_{j}c_{j}}$ is complex 
$\mathcal{K}(l^{2}(\mathbb{N}))$-semicircular in 
($\widetilde{\mathcal{A}}, \widetilde{\varphi})$, for any 
$A_{j}\in \mathcal{K}(l^{2}(\mathbb{N}))$. 
\end{cor}
 Let $L^{2}(\mathcal{F}(\widetilde{H})))$
 be the closure of the algebra of polynomials in 
$\{G(\widetilde{x}), \widetilde{x}\in\widetilde{H^0}\}$ with coeficients in $\mathcal{K}(l^{2}(\mathbb{N}))$ 
\ under the norm $\widetilde{x}\mapsto(||\widetilde{\varphi}(\widetilde{x}^{*}\widetilde{x})||)^{\frac{1}{2}}$.

The following Theorem is the analogue of the classical result( \cite{GJ}, \cite{Jan}, \cite{Si}):
\textit{Let $(V,\langle,\rangle)$ be a real Hilbert space of Gaussian  variables and $\mathcal{A}(V)$ be the algebra of polynomials in elements of $V$. Then there is a unique, up to an isometry, real Hilbert space $H_0$ together with an isometry $\psi:{(V,\langle,\rangle)}\rightarrow
(\{G_1(\eta);\eta\in{H_0}\},\langle{\Omega},\cdot{\Omega}\rangle)$ 
which extends to an isometric isomorphism from $\overline{\mathcal{A}(V)}^{\langle\cdot,\cdot\rangle}$ 
to 
$L^{2}(\Gamma_1,\langle{\Omega},\cdot{\Omega}\rangle)$}.
\begin{thm}\label{thm2}
 Let $(\mathcal{A},\varphi)$ be a $\ast$-non-commutative probability space, $V\subset{A}$ a linear subspace of selfadjoint elements,  and $\mathfrak{A}$ the algebra generated by the polynomials in elements of 
$\mathcal{K}(l^{2}(\mathbb{N}))\otimes{V}$. 
Then, with the above notations, we have that all the elements of $(\mathcal{K}(l^{2}(\mathbb{N}))\otimes{V})_{SA}$ 
are $\mathcal{K}(l^{2}(\mathbb{N}))$-semicircular if and only if there is an unique, up to an isometry, real Hilbert space $\mathcal{H}$ together with an isometry $\psi:(V,\varphi)\rightarrow{(\{G_0(\eta);
 \eta\in {\mathcal{H}}\},\langle{\Omega},\cdot{\Omega}\rangle})$ 
which extends to an isometric isomorphism from $\overline{\mathfrak{A}}^{||\widetilde{\varphi}(\cdot,\cdot)||}$ 
to 
$L^{2}(\mathcal{K}(\widetilde{H}))$.
\end{thm}
\begin{proof}
Suppose we have $(\mathcal{H}, \psi)$ as in the hypothesis. Any  $\widetilde{x}\in{\mathcal{K}(l^{2}(\mathbb{N}))\otimes{V})_{SA}}$ 
can be approximated with elements of the form 
$\sum_{j=1}^{m}{A_j}{x_j}$, with $x_j\in{V}$ and $A_j$ in ${M_n(\mathbb{C})_{SA}}$ 
for some $n\in\mathbb{N}$.

Let $\eta_{j}\in{\mathcal{H}}$ such that 
$x_{j}=\psi^{(-1)}(G_0(\eta_{j}))$, $j=1,\dots,m$. 
From [1], $\{G_0(\eta_j)\}_{j=1}^{m}$ satisfy the 0-Wick Theorem, hence so do $\{{x_j}\}_{j=1}^{m}$, therefore $({x_1}\dots{x_m})$ is a semicircular family  and Thm 1. implies $\widetilde{x}$ is 
$\mathcal{K}(l^{2}(\mathbb{N}))$-semicircular, q.e.d..

Conversely, suppose that all the elements from 
$\mathcal{K}(l^{2}(\mathbb{N}))\otimes{V}_{SA}$
 are $\mathcal{K}(l^{2}(\mathbb{N}))$-semicircular.
 For $x\in{V}$ let $\hat{x}=G_0^{(-1)}(\psi(x))$. If $a,b$ are arbitrary elements form $V$ then $\langle{\hat{a}},{\hat{b}}\rangle
=\langle{\psi(a)}\psi(b)\Omega,\Omega\rangle=\varphi(ab)$, 
therefore the uniqueness is proven.

 For the existence, let $V_0$ be the closure of $V$ under $\parallel{x}\parallel=\sqrt{\varphi(x^2)}$. 
Consider $\mathcal{H}=(V_{0}, \langle,\rangle)$
 with $\langle{x},{y}\rangle=\varphi(x^{*}y)$ 
and $x\rightarrow{\hat{x}}$ the natural identification from $V_0$ to $H_0$. Then $\psi(\cdot)=G_0(\hat{\cdot})$ is an isometry.

 Since any selfadjoint 
$\widetilde{x}\in \mathcal{K}(l^{2}(\mathbb{N}))\otimes{V}$ 
can be approximated wtih elements of the form $\sum_{j=1}^{m}{A}_j{{x_j}}
=\sum_{j=1}^{m}{A}_j{\psi^{-1}(G_0({\widehat{x_j}}))}$ 
for some 
${n\in\mathbb{N}}$, $x_1,\dots,x_m\in{V}$ and $A_1,\dots{A_m}\in({M_n(\mathbb{C})\otimes{V})_{SA}}$, 
from the the 0-Wick theorem we have that :
\begin{eqnarray*}
\langle{\widetilde{\Omega},\psi(\widetilde{x_1})\dots}
& &\mspace{-46mu}{{\psi(\widetilde{x_m})})\widetilde{\Omega}\rangle}=\\
&=&
\langle{\widetilde{\Omega},(\sum_{1\leq{j_1}\leq{p_1}}{A_{1,j_1}}\psi(x_{1,j_1}))
\dots(\sum_{1\leq{j_m}\leq{p_m}}{A_{m,j_m}}\psi(x_{m,j_m}))\widetilde{\Omega}\rangle}\\
&=&
\langle{\widetilde{\Omega},\sum_{\substack{(j_1,\dots,j_m)\\
{1\leq{j_k}\leq{p_k}}}}({A_{1,j_1}}\psi(x_{1,j_1}))
\dots({A_{m,j_m}}\psi(x_{m,j_m}))\widetilde{\Omega}\rangle}\\
&=&
\sum_{\substack{(j_1,\dots,j_m)\\{1\leq{j_k}\leq{p_k}}}}{A_{1,j_1}}
\dots{A_{m,j_m}}\langle{\Omega,G_0(\widehat{x_{1,j_1}})\dots{G_0(\widehat{x_{m,j_m}})}
\Omega\rangle}\\
&=&
\sum_{\substack{(j_1,\dots,j_m)
\\{1\leq{j_k}\leq{p_k}}}}{A_{1,j_1}}\dots{A_{m,j_m}}(\sum_{\gamma\in{NCPP(m)}}
(\prod_{(k,l)\in\gamma}\langle{\Omega,G_0(\widehat{x_{k,j_k}}){G_0(\widehat{x_{l,j_l}})}
\Omega\rangle}))\\
&=&
\sum_{\substack{(j_1,\dots,j_m)\\
{1\leq{j_k}\leq{p_k}}}}{A_{1,j_1}}\dots{A_{m,j_m}}(\sum_{\gamma\in{NCPP(m)}}
(\prod_{(k,l)\in\gamma}\langle{\widehat{x_{k,j_k}},\widehat{x_{l,j_l}}\rangle}))\\
&=&\sum_{\substack{(j_1,\dots,j_m)\\
{1\leq{j_k}\leq{p_k}}}}{A_{1,j_1}}
\dots{A_{m,j_m}}(\sum_{\gamma\in{NCPP(m)}}
(\prod_{(k,l)\in\gamma}\varphi({x_{k,j_k}}{x_{l,j_l}})))\\
&=&
\sum_{\substack{(j_1,\dots,j_m)\\{1\leq{j_k}\leq{p_k}}}}{A_{1,j_1}}
\dots{A_{m,j_m}}\varphi({x_{1,j_1}}\dots{x_{m,j_m}})\\
&=&
\sum_{\substack{(j_1,\dots,j_m)\\{1\leq{j_k}\leq{p_k}}}}\varphi_n({A_{1,j_1}}{x_{1,j_1}}
\dots{A_{m,j_m}}{x_{m,j_m}})\\
&=&
\varphi_n((\sum_{1\leq{j_1}\leq{p_1}}{A_{1,j_1}}{x_{1,j_1}})
\dots(\sum_{1\leq{j_m}\leq{p_m}}{A_{m,j_m}}{x_{m,j_m}}))\\
&=&
\varphi_n(\widetilde{x_1}\dots\widetilde{x_m})
\end{eqnarray*}
\end{proof}
\begin{cor}
 Let $(\mathcal{A},\varphi)$ be a $\ast$-non-commutative probability space, $V\subset{A}$ a linear subspace of selfadjoint elements,  and $\mathfrak{A}$ the algebra generated by the polynomials in elements of 
$\mathcal{K}(l^{2}(\mathbb{N}))\otimes{V}$. 
On 
$\mathcal{K}(l^{2}(\mathbb{N}))\otimes{V}$
 we define the $\mathcal{K}(l^{2}(\mathbb{N})$-valued inner product 
$\langle\widetilde{x},\widetilde{y}\rangle=
\widetilde{\varphi}(\widetilde{x}^*\widetilde{y})$. 
Then, 
$\mathcal{K}(l^{2}(\mathbb{N}))\otimes{V}$
 is a inner-product $\mathcal{K}(l^{2}(\mathbb{N}))$-bimodule of complex $\mathcal{K}(l^{2}(\mathbb{N}))$-semicircular elements  if and only if  there is an unique, up to an isometry, real Hilbert space $\mathcal{H}$ together with an isometry $\psi:(V,\varphi)\rightarrow{(\{G_0(\eta); \eta\in {\mathcal{H}}\},\langle{\Omega},\cdot{\Omega}\rangle})$ 
which extends to an isometric isomorphism from $\overline{\mathfrak{A}}^{||\widetilde{\varphi}(\cdot,\cdot)||}$ 
to 
$L^{2}(\mathcal{F}(\widetilde{H}))$
\end{cor}
\begin{proof}
 Suppose we have $(\mathcal{H}, \psi)$ as in the hypothesis. Then any $\widetilde{x}\in{\mathcal{K}(l^{2}(\mathbb{N}))\otimes{V}}$
 can be approximated in the C*-norm by elements of the form:
\begin{eqnarray*}
\widetilde{x_{\iota}}
&=&
\sum_{j=1}^{m}A_jx_j,\qquad
 \text{for some}\ n\in\mathbb{N}, A_j\in{M_n(\mathbb{C})}, x_j\in{V} \\
&=&
\sum_{j=1}^{m}\Re(A_j)x_j+i\sum_{j=1}^{m}\Im(A_j)x_j,\ 
\Re(A_j), \Im(A_j)\in{M_n(\mathbb{C})_{SA}}, x_j\in{V}\\
&=&
\widetilde{x_1}+i\widetilde{x_2},\qquad 
\widetilde{x_1},\widetilde{x_2}\in({M_n(\mathbb{C})
\otimes{V})_{SA}}\subset(\mathcal{K}(l^{2}(\mathbb{N}))\otimes{V})_{SA}
\end{eqnarray*}

 Consider now 
$B_1, B_2\in{\mathcal{K}(l^{2}(\mathbb{N}))}$ 
such that 
$B_1\widetilde{x_1}+B_2\widetilde{x_2}$
 is selfadjoint. Then
\begin{eqnarray*}
B_1\widetilde{x_1}+B_2\widetilde{x_2}
&=&
\sum_{j=1}^{m}B_{1,j}x_j+\sum_{j=1}^{m}B_{2,j}x_j,\  B_{1,j}, B_{2,j}\in{\mathcal{F}(l^{2}(\mathbb{N}))_{SA}},\  x_j\in{V}\\
&=&
\sum_{j=1}^{m}(B_{1,j}+B_{2,j})\psi^{-1}(G_0((\eta_j))), \quad 
G_0^{-1}(\psi(x_j))=\eta_j\in{\mathcal{H}}\\
&=&
\sum_{j=1}^{m}\psi^{-1}(G(B'_j\eta_j)), \quad 
B'_j\in{\mathcal{K}(l^{2}(\mathbb{N}))_{SA}},\ \eta_j\in\mathcal{H}\\
&=&
\psi^{-1}(G(\widetilde{\eta})),\quad \widetilde{\eta}=\sum_{j=1}^{m}B'_j\eta_j\in{\widetilde{H_n^0}}
\end{eqnarray*}
 and from $\boldsymbol{4.3}$ we have q.e.d..

 To prove that 
$\mathcal{K}(l^{2}(\mathbb{N}))\otimes{V}$ 
is an inner-product bimodule, fix $\{e_j\}_{j\in{J}}$ a orthonormal basis in $\mathcal{H}$. Any $\widetilde{x}\in{\mathcal{K}(l^{2}(\mathbb{N}))\otimes{V}}$ 
can be approximated with as 
$\widetilde{x_{\iota}}=\sum_{j=1}^{m}A_j\psi^{-1}(G(e_j))$, 
with 
$A_j\in{\mathcal{F}(l^{2}(\mathbb{N}))}$, 
therefore
\begin{eqnarray*}
\varphi_n(\widetilde{x}^*\widetilde{x})
&=&
id_{\mathcal{K}(l^{2}(\mathbb{N}))}
\otimes{\varphi}((\sum_{j=1}^{m}A_j\psi^{-1}(G(e_j)))^*
(\sum_{j=1}^{m}A_j\psi^{-1}(G(e_j))))\\
&=&
\sum_{i,j=1}^{m}A_i^*A_j\varphi(\psi^{-1}(G(e_i))\psi^{-1}(G(e_j)))\\
&=&
\sum_{i,j=1}^{m}A_i^*A_j\langle{e_i,e_j}\rangle\\
&=&
\sum_{j=1}^{m}A_j^*A_j\geq{0},\quad
\text{with equality for}\  \widetilde{x}=0\\
\end{eqnarray*}
  
 Conversely, if $\mathcal{K}(l^{2}(\mathbb{N}))\otimes{V}$ is an inner-product bimodule of $\mathcal{K}(l^{2}(\mathbb{N}))$-complex semicircular elements, then  $(\mathcal{K}(l^{2}(\mathbb{N}))\otimes{V})_{SA}$ 
is a real liniar space of 
$\mathcal{K}(l^{2}(\mathbb{N}))$-semicirculars and the conclusion follows from $\ref{thm2}$.
\end{proof}

\end{document}